\documentclass[11pt,a4paper]{article}
\usepackage[english]{babel}
\usepackage{amsmath, amsthm, amsfonts, amssymb, amscd}
\usepackage{comment}
\usepackage[utf8]{inputenc}
\usepackage[pdftex]{graphicx}
\usepackage[matrix,arrow]{xy}
\usepackage{enumerate}

\usepackage{tikz}
\usepackage{tkz-fct}

\DeclareGraphicsExtensions{.bmp,.png,.pdf,.jpg}

\newtheorem{teo}{Theorem}

\newtheorem{cor}[teo]{Corollary}
\newtheorem{prop}[teo]{Proposition}

\DeclareMathOperator{\car}{char}
\DeclareMathOperator{\pr}{pr}
\DeclareMathOperator{\spec}{Spec}
\DeclareMathOperator{\ord}{ord}
\DeclareMathOperator{\mult}{mult}
\DeclareMathOperator{\D}{\mathcal D}
\newcommand{\Oh}{\mathcal{O}}
\newcommand{\pf}{\noindent{\bf Proof\ \ }}
\newcommand{\cqd}{{\hfill $\rule{2mm}{2mm}$}\vspace{3mm}}
\newcommand{\Z}{\mathbb Z}
\newcommand{\N}{\mathbb N}
\newcommand{\tildeOh}{\widetilde{\mathcal{O}}}
\newcommand{\I}{\mathcal{I}}
\newcommand{\J}{\mathcal{J}}
\newcommand{\T}{\mathcal{T}}
\newcommand{\R}{\mathcal{R}}
\newcommand{\C}{\mathcal{C}}
\newcommand{\K}{\mathcal{K}}

\DeclareMathOperator{\RM}{RM}
\DeclareMathOperator{\AM}{AM}
\DeclareMathOperator{\M}{M}

\usepackage{float}

\begin{document}
\title{Tjurina number of a local complete intersection curve}

\author{V. Bayer, \ E. Guzm\'an,  \ A. Hefez \ and \ M. Hernandes}

\date{ }
\maketitle

\noindent{\bf Abstract}
Differently from the Milnor number, no formula relating the Tjurina number of a reducible algebroid curve to invariants of its branches was known. The aim of this work is to provide such a formula for a complete intersection algebroid curve with several branches in terms of invariants of its branches and the maximal points of the set of values of K\"ahler differentials on components of the curve.

\bigskip

\medskip

\bigskip

Let $k$ be an algebraically closed field of characteristic zero and let $C=\spec(\Oh)$ be an algebroid curve, where $\Oh=k[[X_1,\ldots,X_n]]/\mathfrak A$ is a one dimensional reduced $k$-algebra with embedding dimension $n$.

Let 
$$
\Omega(\Oh) = \frac{\Oh^{\oplus n}}{\langle f_{X_1}e_1+\cdots+f_{X_n}e_n; \ f\in \mathfrak A\rangle},
$$
where $e_1,\ldots,e_n$ is the canonical basis of $\Oh^{\oplus n}$ and $f_{X_i}=\frac{\partial f}{\partial X_i}$, $i=1,\ldots,n$, be the $\Oh$-module of K\"ahler differentials on $C$. 

We will denote the image of $X_i$ in $\Oh$ by $x_i$ that of $e_i$ in $\Omega(\Oh)$ by $dx_i$. For $g\in \Oh$ we consider the universal differentiation $dg=g_{X_1}dx_1+\cdots+g_{X_n}dx_n\in \Omega(\Oh)$. 

The jacobian ideal $\mathfrak J$ of $\Oh$ is the smallest nonzero Fitting ideal of $\Omega(\Oh)$, that is, the nonzero ideal generated by the minors of greatest possible size of the Jacobian matrix 
\[
Jac(f_1,\ldots,f_m)= \begin{bmatrix} f_{1X_1} & \cdots &  f_{1X_n}\\
                                  \vdots&           & \vdots\\
                                   f_{mX_1} & \cdots & f_{mX_n}
	\end{bmatrix},
	\]				
	where $f_1,\ldots, f_m$ are generators of the ideal $\mathfrak A$. 
													
Let $\wp_1,\ldots,\wp_r$ be the minimal primes of $\Oh$. We set $I=\{1,\ldots,r\}$. For $J=\{j_1<j_2<\cdots<j_s\}\subset I$, define $\Oh_J=\Oh/\cap_{j\in J} \wp_{j}$ and denote by $\pi_J\colon \Oh \to \Oh_J$ the natural projection. When $J=\{i\}$, we denote $\Oh_J$ simply by $\Oh_i$ and $\pi_J$ by $\pi_i$. Since $\Oh$ is reduced, we have that $\Oh_I=\Oh$. The maps $\pi_i$ determine an $\Oh$-modules monomorphism $\pi\colon \Oh \hookrightarrow \prod_{i=1}^r\Oh_i$, which induces an $\Oh$-modules homomorphism $\pi^*\colon \Omega(\Oh) \to \prod_{i=1}^r\Omega(\Oh_i)$.

Let $\widetilde{\Oh_J}$ be the integral closure of $\Oh_J$ in its total ring of fractions $\K_J$. It is well known that 
$
\widetilde{\Oh_J}\simeq \prod_{i=1}^s \widetilde{\Oh_{j_i}} \simeq \prod_{i=1}^s k[[t_{j_i}]]$,
so that $\K_J \simeq \prod_{i=1}^s k((t_{j_i}))$. Define $\K=\prod_{i=1}^r k((t_{i}))$.

Denote by $\varphi\colon \Oh \stackrel{\pi}{\hookrightarrow} \prod_{i=1}^r\Oh_i \hookrightarrow \K$ the natural injective $\Oh$-modules homomorphism and by $\varphi^*$ the induced  $\Oh$-modules homomorphism
$$
\begin{array}{rccccc}
\varphi^* \colon \qquad  \Omega(\Oh) & \stackrel{\pi^*}{\longrightarrow} &   \prod_{i=1}^r\Omega(\Oh_i)& \longrightarrow & \K \\
                  \sum_{l=1}^ng_ldx_l &\mapsto &  \pi^*(\sum_{l=1}^ng_ldx_l) & \mapsto& (\ldots,\sum_{l=1}^n\varphi_i(g_l)\varphi_i(x_l)'t_i, \ldots),		
					\end{array}
$$
where the dash means derivative with respect to $t_i$.

We define the value map $v \colon \K \to (\Z\cup \{\infty\})^r$, where $v_i=\ord_{t_i}$ and where $(\Z\cup \{\infty\})^r$ is endowed with the product order. This allows us, by composition with $\varphi$ and $\varphi^*$, to define value maps on $\Oh$ and on $\Omega(\Oh)$, respectively.

A fractional ideal $\I$ of $\Oh$ is an $\Oh$-submodule of $\K$ such that there exists a regular element $d\in \Oh$ such that $d\;\I\subset \Oh$. We define the value set of $\I$ as being $E=v(\I)\subset (\Z\cup \{\infty\})^r$. In particular, the value set $\Gamma(\Oh)=v(\Oh)$ will be called the semigroup of values of $\Oh$.

Properties of value sets of the form $E=\Gamma(\Oh)$ were studied in \cite{D87}, but it is not difficult to extend them to any value set $E$. We list below some of their remarkable properties. 

A value set $E$ has an infimum $\inf(E)$ (for the product order) and a conductor $c(E)$ defined by
$ c(E) = \inf\{\beta\in \Z^r; \; \forall \alpha\in \Z^r, \alpha\geq \beta \Rightarrow \alpha\in E\}$.


The {\em fiber with respect to a subset} $J\subset I=\{1,\ldots,r\}$ of an element $\alpha\in E$ is the set
$ F_J(E,\alpha)=\{\beta\in E; \beta_j=\alpha_j, \ \text{if} \ j\in J, \ \text{and} \ \beta_i>\alpha_i, \ \text{if} \ i\notin J\}$.
The {\em fiber} of $\alpha$ in $E$ is the set $F(E,\alpha)=\bigcup_{i=1}^rF_{\{i\}}(E,\alpha)$.

We say that $\alpha$ is a {\em maximal point} of $E$ if $F(E,\alpha)=\emptyset$. This means that there is no element in $E$ with one coordinate equal to the corresponding coordinate of $\alpha$ and the others bigger.

From the fact that $E$ has an infimum $\inf(E)$ and a conductor $c(E)$, one has immediately that all maximal points of $E$ are in the limited region $\{x\in \Z^r; \ \inf(E) \leq x < c(E), \ \ i=1,\ldots,r\}$, so $E$ has finitely many maximal points.

A maximal point $\alpha$ of $E$ is called an {\em absolute maximal} if $F_J(E,\alpha)=\emptyset$ for every $J\subset I$, $J\neq I$; and a {\em relative maximal} if $F_J(E,\alpha)\neq\emptyset$, for every $J\subset I$ with $\#J\geq2$.

When $r=2$, maximal, relative maximal and absolute maximal coincide. For $r=3$ there are only relative maximals and absolute maximals. 

We denote by $M(E)$, $RM(E)$ and $AM(E)$ the sets of maximals, of relative maximals and absolute maximals of the set $E$, respectively.

In \cite[Theorem 1.5]{D87} it is proved that the finite set $RM(E)$ and $E_i=\Gamma(\Oh_i)$, $i=1,\ldots,r$, determine $E=\Gamma(\Oh)$ in a combinatorial sense. However, it is not difficult to extend this property for a general value set $E$.

We will denote by $\T$ and by $\T_i$ the torsion modules of $\Omega(\Oh)$ and of $\Omega(\Oh_i)$, respectively.
One has the following:

\begin{prop}\label{irred} Suppose that $C$ is irreducible, that is, $\Oh$ is a domain, then $\ker(\varphi^*)=\T$.
\end{prop}
\pf See \cite[Proposition 1, Sect. 7.1]{HH}.
\cqd

\begin{prop}\label{inclusion} Suppose that $\Oh$ is reduced, then $\pi^*(\T) \subset \T_1\times \cdots \times \T_r$.
\end{prop}
\pf
Let $\omega\in \T$, then there exists an element $h\notin Z(\Oh)=\bigcup_{i=1}^r \wp_i$ such that $h\omega=0$. Since $h\notin Z(\Oh)$, it follows that $\pi_j(h)\neq 0$, for all $j=1,\ldots,r$, hence a nonzero divisor in $\Oh_j$. Since $0=\pi_j^*(h\omega)=\pi_j(h) \pi_j^*(\omega)$, this implies that $\pi_j^*(\omega)\in \T_j$, proving the result.
\cqd

We will make an essential use of the following theorem. 

\begin{teo}[Piene \cite{Pi}] \label{piene} If $\Oh$ is the ring of an algebroid complete intersection curve, then
\[
\mathfrak J\tildeOh=\mathcal{C}\R,
\]
where $\mathfrak J$ and $\mathcal{C}$ are, the Jacobian ideal of $\Oh$ and the conductor ideal $\Oh:\tildeOh$, respectively, while $\R$ is the ramification ideal generated in $\tildeOh$ by the elements $ (\varphi_1(x_i)',\ldots, \varphi_r(x_i)')$,\;  $i=1,\ldots,n$. 
\end{teo}

If we denote $c(\Gamma(\Oh))$ by $c$, then one has that $\C=\{g\in \K; \; v(g)\geq c\}$.\smallskip

From now on we will assume that the curve $C$ is a complete intersection, this means that $C=\spec k[[X_1,\ldots,X_n]]/\langle f_1,\ldots,f_{n-1}\rangle$, for some $f_1,\ldots,f_{n-1}\in k[[X_1,\ldots,X_n]]$.\medskip

By choosing properly the coordinates in $k[[X_1,\ldots,X_n]]$, we may assume that the hyperplane $X_1$ does not contain any line of the tangent cone of the curve $C$, so $(\varphi_1(x_1),\ldots,\varphi_r(x_1))$ is not in the set $Z(\tildeOh)$ of zero divisors of $\tildeOh$; hence the maximal ideal $\mathcal M$ of $\Oh$ is a fractional ideal of $\tildeOh$ and 
$$v(x_1)=(\mult(C_1),\ldots,\mult(C_r))= \inf v(\mathcal M).$$

Since $\car(k)=0$, it follows that $\R$ is a regular fractional ideal of $\Oh$, that is, a fractional ideal that contains a regular element, and 
$$\inf v(\R)= v(\varphi_1(x_1)',\ldots,\varphi_r(x_1)')=v(x_1)-e, \ \text{where} \ e=(1,\ldots,1).$$

In this situation, one has that $\mathfrak J=\langle |M_1|,\ldots, |M_n| \rangle$, where $|M_i|$ is the determinant of the matrix $M_i$ obtained by deleting the $i$-th column of the matrix $Jac(f_1,\ldots,f_{n-1})$. From Theorem \ref{piene}, after reordering the series $f_1,\ldots,f_{n-1}$, we may assume that  
\begin{equation}\label{eq1}
v(|M_1|)=\inf v(\mathfrak J\tildeOh) =\inf v(\mathcal{C}\R)=c+ \inf v(\R)=c+v(x_1)-e.
\end{equation}	
This implies that $|M_1|$ is not a zero divisor in $\K$.		\medskip																												

Writing the relations $df_i=\sum_{j=1}^nf_{iX_j}dx_j=0, \ \ i=1,\ldots,n-1$, in $\Omega(\Oh)$ in matricial form 
\[ 
M_1\begin{bmatrix} dx_2 \\ \vdots \\ dx_n \end{bmatrix} =-\begin{bmatrix}f_{1X_1}\\ \vdots \\ f_{n-1X_1} \end{bmatrix}\;dx_1,
\] 
we get, by Cramer's Rule,
\begin{equation}\label{eq2}
|M_1|\;dx_i=\sigma(i) |M_i|\;dx_1, \ \text{where} \ \sigma(i)=(-1)^i
\end{equation}

From Equations (\ref{eq1}) and (\ref{eq2}) one gets that, for all $i=1,\ldots,r$,
\begin{equation}\label{eq3}
v(|M_i|)=c+v(x_i)-e.
\end{equation}

\begin{teo} \label{torsion} For a complete intersection $\Oh$ one has that $$(\pi^*)^{-1}(\T_1\times \cdots\times \T_r)=\T.$$
\end{teo}
\pf
From Proposition \ref{inclusion} one has that $\T \subset (\pi^*)^{-1}(\T_1\times \cdots \times \T_r)$.
 For the other inclusion, let $\omega=\sum_{l=1}^ng_ldx_l \in (\pi^*)^{-1}(\T_1\times \cdots\times \T_r)$. Hence $\pi_i^*(\omega)\in \T_i$, therefore $\varphi_i^*(\omega)=0$, for all $i=1,\ldots,r$ (cf. Proposition \ref{irred}). 

Consider $|M_1|\;\omega$, which in view of relations (\ref{eq2}) may be written as
\begin{equation} \label{omegatorsion}
|M_1|\;\omega=\sum_{l=1}^n g_l|M_1|\;dx_l=(\sum_{l=1}^n \sigma(l) g_l|M_l|)\;dx_1.
\end{equation}

Applying $\varphi_i^*$ to both sides of (\ref{omegatorsion}), one gets, for all $i=1,\ldots,r$,
\[
0 =  \varphi_i^*(|M_1|\;\omega) = \varphi_i\big(\sum_{l=1}^n \sigma(l) g_l|M_l|\big) \varphi_i(x_1)'\;t_i.
\]

Since $\varphi_i(x_1)'\;t_i\neq 0$, for all $i$, then $\varphi_i (\sum_{l=1}^n \sigma(l) g_l|M_l|)=0$, for all $i$, which implies that  $\sum_{l=1}^n \sigma(l) g_l|M_l| \in \bigcap_{i} \wp_i=(0)$. So, from (\ref{omegatorsion}), we have $|M_1|\;\omega =0$ and since $|M_1|\notin Z(\Oh)$, it follows that $\omega\in \T$.
\cqd

\begin{cor} Let $\Oh$ be the ring of an algebroid complete intersection curve, then one has $\omega=\sum_{l=1}^ng_ldx_l \in \T$ if and only if \[
\dfrac{df_1\wedge \cdots \wedge df_{n-1}\wedge \omega}{dx_1\wedge\cdots\wedge dx_n} =\det\begin{bmatrix} f_{1X_1} & \cdots &  f_{1X_n}\\
                                  \vdots&           & \vdots\\
                                   f_{n-1X_1} & \cdots & f_{n-1X_n}\\
																	g_1& \cdots &g_n
	\end{bmatrix}=0.
	\]				
\end{cor}

\pf We have seen in the course of the proof of Theorem \ref{torsion} that $\omega\in \mathcal T$ if and only if $|M_1|\;\omega=0$, and this is equivalent, in view of (\ref{omegatorsion}), to the condition $\sum_{l=1}^n \sigma(l) g_l|M_l|=0$, which in turn is equivalent to the vanishing of the determinant in the statement of the corollary.
\cqd

From Proposition \ref{irred} and Theorem \ref{torsion} it follows that, if $C$ is a complete intersection, then we have inclusions
\[
\Omega(\Oh)/\T \hookrightarrow \prod_{i=1}^r \Omega(\Oh_i)/\T_i \hookrightarrow \prod_{i=1}^r \Omega(\tildeOh_i)  \simeq (t_1,\ldots,t_r)\tildeOh \hookrightarrow \K,
\]
which allow us to view $\Omega(\Oh)/\T$ and  $\Omega(\Oh_1)/\T_1 \times \cdots \times \Omega(\Oh_r)/\T_r$ as fractional ideals of $\Oh$, one contained in the other.

Now, it is easy to verify that the  map
\[
\begin{array}{rccc}
\D \colon &\Oh & \to &\Oh\\
&g& \mapsto & \dfrac{df_1\wedge \cdots \wedge df_{n-1}\wedge dg}{dx_1\wedge\cdots\wedge dx_n}
\end{array}
\]
is a differential operator on $\Oh$, that may be written as
\[
\D(g)=\det \begin{bmatrix} f_{1X_1} & \cdots &  f_{1X_n}\\
                                  \vdots&           & \vdots\\
                                   f_{n-1X_1} & \cdots & f_{n-1X_n}\\
																	g_{X_1}& \cdots &g_{X_n}
	\end{bmatrix}.
	\]
	
This operator and the universal differentiation  on $\Oh$ may be extended in a standard way to the total ring of fractions $\mathcal K$ of $\Oh$. 

\begin{prop} For any $g\in \Oh$, and any $i=1,\ldots,r$, one has 
\[
v_i(\D(g))\geq v_i(g)+c_i-1,
\]
with equality if $\pi_i(g)$ is not a unit in $\Oh_i$.
\end{prop}
\pf Since $df_j=0$, for all $j=1,\ldots,n-1$,  from the relations $\varphi_i^*(df_j)=0$, $i=1,\ldots,r$, and the chain rule applied to $\varphi_i(g)$, one gets the following equalities:
\[
\begin{bmatrix} \varphi_i(f_{1X_1}) & \cdots &  \varphi_i(f_{1X_n})\\
                                  \vdots&           & \vdots\\
                                   \varphi_i(f_{n-1X_1}) & \cdots & \varphi_i(f_{n-1X_n})\\
																	\varphi_i(g_{X_1})& \cdots &\varphi_i(g_{X_n})
	\end{bmatrix} \begin{bmatrix} \varphi_i(x_1)'t_i \\ \vdots \\ \varphi_i(x_{n-1})'t_i \\ \varphi_i(x_n)'t_i \end{bmatrix}=
	\begin{bmatrix} 0 \\ \vdots \\ 0\\ \varphi_i(g)'t_i \end{bmatrix}.
	\]
	
	Again, from Cramer's Rule, for $1\leq i\leq r$ and $1\leq j\leq n$, one gets that
	\[
	\varphi_i(\D(g))\varphi_i(x_j)'t_i= \sigma(n+j) \varphi_i(|M_j|)\varphi_i(g)'t_i.
	\]
	
	From this, in view of Equation (\ref{eq3}), we get that
	\[
	v_i(\varphi_i(\D(g)))=c_i+v_i(\varphi_i(g)')\geq c_i+v_i(g)-1,
	\]
	with equality if $\pi_i(g)$ is not a unit in $\Oh_i$.
	\cqd
	
	\begin{prop} One has that $\D(\tildeOh)\subset \tildeOh$.
	\end{prop}
	\pf
	If $\frac{g}{h}\in \tildeOh$, with $g,h\in \Oh$, then $v(g)\geq v(h)$. Since  
	\[
	\D\left(\frac{g}{h}\right)=\frac{h\D(g)-g\D(h)}{h^2},
	\]
	\[
	v\left(\frac{\D(g)}{h}\right) \geq v(g)-v(h)+c-e, 
	\]
	and
	\[v\left(\frac{g\D(h)}{h^2}\right) \geq v(g)-v(h)+c-e,
	\]
	one has that 
\[
v\left(\D\left(\frac{g}{h}\right)\right) \geq v(g)-v(h)+c-e \geq c-e \geq (0,\ldots,0).
\]															
\cqd
																
Using the determinantal definition of $\D(u)$ and the relations (\ref{eq2}) one gets
\[
|M_1|\;du= \sigma(n-1) \D(u)\;dx_1.
\]

So, one has that $\D(u)$ is not a zero divisor, because, otherwise, $|M_1|\;du$ would be a torsion differential, that implies $du$ is torsion, a contradiction, since in caracteristic zero, no exact non-zero differential is torsion. So, 
\[
\xi=\frac{\varphi^*(du)}{\varphi(\D(u))}= \sigma(n-1) \frac{1}{\varphi(|M_1|)}\varphi^*(dx_1) \in \K,
\]
is well defined for $u\not\in k$ and independent from  $u$; furthermore, $v(\xi)=-c+e$. 

Consider now the $\Oh$-modules homomorphism multiplication by $\xi$:
\[ \begin{array}{rccl}
m_\xi\colon & \Oh & \to     & \K .\\
               &  h  & \mapsto & \varphi(h) \xi 
									\end{array}
									\]

Recalling that $\mathfrak J$ and $\C$ are respectively the jacobian and conductor ideals of $\Oh$, we have the following theorem:

\begin{teo} \label{basic} The following assertions are true
\begin{enumerate}[\rm i)]
\item $m_\xi$ is injective;
\item $m_\xi(\mathcal{C})=(t_1,\ldots,t_r)\tildeOh=\varphi^*(\Omega(\tilde\Oh))$,
\item $m_\xi(\mathfrak J)=\varphi^*(\Omega(\Oh))$.
\end{enumerate}
\end{teo}
\pf
(i) This is so, because the homomorphism $\varphi\colon \Oh \to \tildeOh$ is injective, $\xi\neq 0$  and $\K \simeq \prod_{i=1}^rk((t_i))$ is free.
													
\noindent (ii) Since $v(\xi)=-c+e$, it follows that for $h\in \mathcal{C}$, $m_\xi(h)=\varphi(h)\xi$ is such that $v(m_\xi(h))=v(h)+v(\xi)\geq e$, hence $m_\xi(h)\in (t_1,\ldots,t_r)\tildeOh$.

Conversely, let $\rho= (\rho_1,\ldots,\rho_r) \in (t_1,\ldots,t_r)\tildeOh$. Let $h=(h_1, \ldots, h_r)$ , where $h_i=\rho_i/\xi_i$, then $h\in \mathcal{C}$ and $m_\xi(h)=\rho$. 

\noindent (iii) Notice that any element in $\mathfrak J$ may be written in the form
\[
h=\sum_{l=1}^{n} \sigma(l)g_l|M_l|, \ \ \text{for} \ g_l\in \Oh,
\]
hence
\[ 
\begin{array}{rcl}
\left(m_\xi(h)\right)_i&=& \sigma(n-1) \sum_{l=1}^n \dfrac{\sigma(l)\varphi_i(g_l|M_l|)}{\varphi_i(|M_1|)}\varphi_i^*(dx_1)\\ \\
           &=& \sigma(n-1)  \sum_{l=1}^n\varphi_i(g_l)\varphi_i^*(dx_l)\\ \\
					&=& \sigma(n-1)\varphi_i^*(\sum_{l=1}^ng_ldx_l),
					\end{array}
\]
therefore,
\[
m_\xi(h)=\sigma(n-1)\varphi^*(\sum_{l=1}^ng_ldx_l)\in \varphi^*(\Omega(\Oh)).
\]

Conversely, let $\varphi^*(\sum_{l=1}^ng_ldx_l)\in \varphi^*(\Omega(\Oh))$, then defining
\[
h=\sigma(n-1)\sum_{l=1}^{n} \sigma(l)g_l|M_l|\in \mathfrak J,
\]
we get $m_\xi(h)= \varphi^*(\sum_{l=1}^ng_ldx_l)$, concluding the proof. \cqd

\begin{cor}[Pol \cite{Po}, Proposition 3.31] For a complete intersection $\Oh$, one has the following $\Oh$-module isomorphisms
\[
\mathfrak J \simeq \varphi^*(\Omega(\Oh))\simeq \Omega(\Oh)/\T,
\]
and 
\[
v(\mathfrak J)=v(\Omega(\Oh))+c-e.\]
\end{cor}
\pf
This is just Theorem \ref{basic} (iii) and a computation of values.
\cqd

The Tjurina number of $\Oh$ is defined as
\[
\tau(\Oh) = \ell_{\Oh}\left( \dfrac{\Oh}{\mathfrak J} \right),   
\]
where $\ell_\Oh$ stands for the length as $\Oh$-module.

\begin{cor} \label{cor10} Let $\Oh$ be the ring of a complete intersection curve, then one has
\[
\tau(\Oh)= \ell_{\Oh}\left( \dfrac{\Oh}{\mathcal{C}} \right) +\ell_{\Oh}\left( \dfrac{\varphi^*(\Omega(\tildeOh))}{\varphi^*(\Omega(\Oh))}\right).
\]
\end{cor}
\pf By computing values, it is easy to check that $\mathfrak J\subset \mathcal{C}$, so
\[
\tau(\Oh)=\ell_{\Oh}\left( \dfrac{\Oh}{\mathfrak J} \right) = \ell_{\Oh}\left( \dfrac{\Oh}{\mathcal{C}} \right) +\ell_{\Oh}\left(\dfrac{\mathcal{C}}{\mathfrak J} \right).
\]

Now, from Theorem \ref{basic} one gets
\[
\ell_{\Oh} \left( \dfrac{\mathcal{C}}{\mathfrak J}\right)= \ell_{\Oh} \left(\dfrac{m_\xi(\mathcal{C})}{m_\xi(\mathfrak J)}\right)=\ell_{\Oh} \left(\dfrac{\varphi^*(\Omega(\tildeOh))}{\varphi^*(\Omega(\Oh))}\right),
\]
which proves the result.
\cqd

Let us denote by $\Gamma_i$ the semigroup of values of the branch $\Oh_i$, by $\Lambda_i$ the set of values of the fractional ideal $\Omega_i/{\mathcal T}_i$ of $\Oh_i$ and by $\Lambda$ the set of values of the fractional ideal $\Omega/{\mathcal T}$ of $\Oh$. Recall the definition of the singularity degree $\delta=\ell_\Oh\;(\tildeOh/\Oh)$ of $\Oh$. We also define $r_i=\# \Lambda_i\setminus \Gamma_i$. This integer is equal to the number of nonexact linearly independent differentials modulo the exact differentials of the branch $\Oh_i$ that  may be made more explicit when $\Oh_i$ is also a complete intersection: $r_i=c_i-\tau_i$, where $c_i=2\delta_i$ is the conductor of $\Gamma_i$ and $\tau_i$ is the Tjurina number of $\Oh_i$ (cf. \cite[Note on page 783]{Za}.

\begin{cor} Let $\Oh$ be the ring of a complete intersection curve, then one has
\[
\tau(\Oh)= \delta+ \sum_{i=1}^r (\delta_i-r_i) + \ell_{\Oh}\left(\frac{\prod_{i=1}^{r} \frac{\Omega_i}{\mathcal T_i}}{\frac{\Omega(\Oh)}{\mathcal T}}\right).
\]
\end{cor}
\pf  Since $\Oh$ is Gorenstein, it follows that 
\[
 \ell_{\Oh} \left( \dfrac{\Oh}{\mathcal{C}}\right)= \ell_{\Oh} \left( \dfrac{\tildeOh}{\Oh}\right)=\delta.\]

On the other hand, from the chain of monomorphisms
\[
\varphi^*(\Omega(\Oh))\simeq \frac{\Omega(\Oh)}{\T} \ \longrightarrow \ \prod_{i=1}^{r} \frac{\Omega_i}{\mathcal T_i} \ \longrightarrow \ \varphi^*(\Omega(\tildeOh)),
\]
we get that
\[
\ell_{\Oh} \left(\dfrac{\varphi^*(\Omega(\tildeOh))}{\varphi^*(\Omega(\Oh))}\right)= \ell_{\Oh} \left(\dfrac{\varphi^*(\Omega(\tildeOh))}{\prod_{i=1}^{r} \frac{\Omega_i}{\mathcal T_i}}\right)+\ell_{\Oh}\left(\frac{\prod_{i=1}^{r} \frac{\Omega_i}{\mathcal T_i}}{\frac{\Omega(\Oh)}{\mathcal T}}\right).
\]

Since
\[
\ell_{\Oh} \left(\dfrac{\varphi^*(\Omega(\tildeOh))}{\prod_{i=1}^{r} \frac{\Omega_i}{\mathcal T_i}}\right)=\sum_{i=1}^r\ell_{\Oh_i}\left( \frac{\varphi_i^*(\Omega(\tildeOh_i))}{\frac{\Omega_i}{\tilde T_i}}\right),
\]
and
\[\begin{array}{lcl}
\ell_{\Oh_i}\left( \frac{\varphi_i^*(\Omega(\tildeOh_i))}{\frac{\Omega_i}{\mathcal T_i}}\right)=\#\left( v_i(\varphi_i^*(\Omega(\tildeOh_i))) \setminus v_i(\frac{\Omega_i}{\mathcal T_i})\right) &=&\# \N^*\setminus \Lambda_i\\ 
&=&\#\N \setminus \Gamma_i -\# \Lambda_i\setminus \Gamma_i\\
&=& \delta_i - r_i,\end{array}
\]
the proof follows in view of the formula given in Corollary \ref{cor10}.

\cqd

So, to compute $\tau(\Oh)$, it remains only to compute $\ell_{\Oh}\left(\frac{\prod_{i=1}^{r} \frac{\Omega_i}{\mathcal T_i}}{\frac{\Omega(\Oh)}{\mathcal T}}\right)$, which we will do next.\medskip

Let us fix some notation. Let $\I$ be a fractional ideal with value set $E$ and  $\alpha^0=\inf(E)$. For $\gamma\in \Z^r$, we define $\I(\gamma)=\{x\in \I; \; v(x)\geq \gamma\}$. For $J=\{j_1< \cdots < j_s\}\subset I=\{1,\ldots,r\}$, we also denote by $\pi_J\colon \K \to \K_J$ the natural projection and by $\pr_J\colon \Z^r \to \Z^s$ the projection $(\alpha_1,\ldots,\alpha_r)\mapsto (\alpha_{j_1},\ldots,\alpha_{j_s})$. It is true that the value set of $\pi_J(\I)$ is $E_J=\pr_J(E)$. For $1\leq m \leq r$, Let us put $[1,m]=\{1,\ldots, m\}=[1,m+1)$. We also denote by $\pr_*(\alpha_1,\ldots,\alpha_s)=\alpha_s$ and by $\mathcal G(E_i)=\{z\in \mathbb{Z}\setminus E_i; z\geq \inf(E_i)\}$ the set of gaps of $E_i=\pr_{\{i\}}(E)$.

If $\J \subset \I$ are fractional ideals of $\Oh$, and for $\gamma\in \Z^r$ sufficiently large (in the order product of $\Z^r$), we have that
\[
\ell_\Oh\left(\frac{\I}{\J}\right)= \ell_\Oh\left(\frac{\I}{\I(\gamma)}\right)-\ell_\Oh\left(\frac{\J}{\J(\gamma)}\right).
\]

In particular, for such a $\gamma$, we have that
\begin{equation} \label{diff}
 \ell_{\Oh}\left(\frac{\prod_{i=1}^{r} \frac{\Omega_i}{\mathcal T_i}}{\frac{\Omega(\Oh)}{\mathcal T}}\right)= \ell_{\Oh}\left( 
\frac{\prod_{i=1}^{r} \frac{\Omega_i}{\mathcal T_i}}{\prod_{i=1}^{r} \frac{\Omega_i}{\mathcal T_i}(\gamma)}\right)-\ell_{\Oh}\left(\frac{\frac{\Omega(\Oh)}{\mathcal T}}{\frac{\Omega(\Oh)}{\mathcal T}(\gamma)}\right).
\end{equation}\medskip

Let us define
\[
\Theta_1=0, \quad \text{and} \quad  \Theta_i= \# \bigcup_{\{ i\} \subsetneq J\subseteq [1,i]}\pr_*(\RM(\Lambda_J)), \quad 2\leq i\leq r,
\]
and recall the following result, adapted to our situation:

\begin{teo}[Guzm\'an-Hefez \cite{GH}, Theorem 10] \label{teogeral1} Let $\I$ be a fractional ideal of the ring $\Oh$ of an algebroid curve with values set $E$. Suppose that $\Oh$ has $r$ minimal primes. If $\gamma$ is sufficiently large, then
\[ 
\ell\left(\dfrac{\I }{\I(\gamma)}\right)=\ell\left(\dfrac{\pi_{[1,r)}(\I)}{\pi_{[1,r)}(\I)(pr_{[1,r)}(\gamma))}\right)+(\gamma_r-\alpha^0_r)-\#\mathcal G (E_r)- \Theta_r.\] 
\end{teo}

Now, since $\prod_{i=1}^{r} \Lambda_i$ has no maximal points at all, we have
\begin{equation}\label{ellprod}
\ell_{\Oh}\left( \frac{\prod_{i=1}^{r} \frac{\Omega_i}{\mathcal T_i}}{\prod_{i=1}^{r} \frac{\Omega_i}{\mathcal T_i}(\gamma)}\right)=\sum_{i=1}^r (\gamma_i-\alpha^0_i -\# {\mathcal G}(\Lambda_i)).
\end{equation}

Applying Theorem \ref{teogeral1} and Formulas (\ref{diff}) and (\ref{ellprod}) to $\I=\Omega(\Oh)/\mathcal T$ and $\I_i=\Omega_i/\mathcal T_i$, we get
\[ 
\ell_{\Oh}\left( \frac{\prod_{i=1}^{r} \frac{\Omega_i}{\mathcal T_i}}{\frac{\Omega(\Oh)}{\mathcal T}}\right) =
\sum_{i=1}^r (\gamma_i-\alpha^0_i -\# {\mathcal G}(\Lambda_i) )-  \sum_{i=1}^r (\gamma_i-\alpha^0_i -\# {\mathcal G}(\Lambda_i) - \Theta_i)= \sum_{i=1}^r \Theta_i
 \]

So, we proved the following result:

\begin{teo} Let $\Oh$ be a complete intersection curve, then
\[
\tau(\Oh)= \delta+ \sum_{i=1}^r (\delta_i-r_i +\Theta_i).
\]
\end{teo}

For $r=2$ and $r=3$, from \cite[Proposition 7]{GH} and \cite[Theorem 18]{GH}, respectively, one gets the more explicit formulas: 
\[
\begin{array}{rl}
\tau(\Oh)= & \delta+ \delta_1+\delta_2-(r_1+r_2) +\# \M(\Lambda),\\ \\
\tau(\Oh)=&\delta+ \delta_1+\delta_2+\delta_3-(r_1+r_2+r_3) +\sum_{1\leq i<j\leq 3} \# \M(\Lambda_{\{i,j\}})  + \\
          &\# \RM(\Lambda)- \# \AM(\Lambda).
					\end{array}
\]

\begin{cor} Suppose that $\Oh$ and each $\Oh_i$ are complete intersections (this is the case for plane curves), then
\[
\tau(\Oh)= \sum_{i=1}^r ( \tau_i + \frac{1}{2} I_i+\Theta_i) .
\]
\end{cor}
\pf This is so, because $2\delta_i- r_i=\tau_i$ and due to the formula (cf. \cite{H})
\[
\delta=\sum_{i=1}^r (\delta_i + \frac{1}{2}  I_i),   \ \ \text{where} \ \ I_i=\dim_k \frac{\Oh}{\wp_i +\cap_{l\neq i} \wp_l}.
\]
\cqd

In particular, for $r=2$, our formula specializes as
\[
\tau(\Oh)=\tau_1+\tau_2+\frac{1}{2}(I_1+I_2)+\# \M(\Lambda);
\]
and for $r=3$, as
\begin{equation}\label{last}
\begin{array}{rl}
\tau(\Oh)=&\tau_1+\tau_2+\tau_3 +\frac{1}{2}(I_1+I_2+I_3) + \sum_{1\leq i<j\leq 3} \# \M(\Lambda_{\{i,j\}})  + \\
         &  \# \RM(\Lambda)- \# \AM(\Lambda).
\end{array}\end{equation}\medskip

\begin{cor}[Additive formula] Let $f=f_1\cdots f_r\in k[[X,Y]]$ be the equation of an algebroid reduced plane curve. Let $J,K\subset I=\{1,\ldots,r\}$ be a partition of $I$. Then one has 
\[
\tau(\Oh)=\tau(\Oh_I)=\tau(\Oh_J)+\tau(\Oh_K)+\sum_{j\in J, \kappa\in K}I(f_j,f_\kappa) + \Theta(\Oh_I)-(\Theta(\Oh_J)+\Theta(\Oh_K)),
\]
where $I(f_j,f_\kappa)$ stands for the intersection multiplicity of $f_j$ and $f_\kappa$, and for any $J=\{j_1<\cdots<j_s\} \subset I$,
\[
\Theta(\Oh_J)=\sum_{i=2}^s \# \bigcup_{\{j_i\}\subsetneq L \subset \{j_1,\ldots,j_i\}}\pr_*(\RM(\Lambda_L)).
\]
\end{cor}
\bigskip

\begin{center} {\bf Example} \end{center}

\medskip

In what follows we give an example where our formula is applied. 
\medskip

Let  $\Oh=k[[X,Y]]/\langle f\rangle$, where
$$f(X,Y)=(Y^3-X^7)(Y^3-3X^5Y-X^7-X^8)(Y^4-2X^5Y^2-4X^7Y-X^9+X^{10}).$$

By using the algorithm decribed in \cite{HH2}, we get 
\[ \begin{array}{ll}
\Lambda_1=v(\Omega_1/\T_1)=\{3,6,7,9,10,12,\ldots\},&   c(\Lambda_1)=12,\\ \\
\Lambda_2=v(\Omega_2/\T_2)=\{3,6,7,9,10,12,\ldots\},&   c(\Lambda_2)=12,\\ \\
\Lambda_3=v(\Omega_3/\T_3)=\{4,8,9,12,13,14,16,\ldots\}, & c(\Lambda_3)=16.
\end{array}
\]

For the computations below, we used a script implemented in MAPLE by the last author.\medskip

The set of maximal points of $\Lambda_{\{1,2\}}$, whose conductor is $(20,20)$, is 
\[\begin{array}{rl} 
\M(\Lambda_{\{1,2\}})=&\{(3,3),(6,6),(7,7),(9,9),(10,10),(12,12),(13,13), (15,15),\\
                      &(16,16),(19,19)\}.
																					\end{array}
\]
Notice that  $\#\M(\Lambda_{\{1,2\}})=10$.\medskip

The set of maximal points of $\Lambda_{\{1,3\}}$, whose conductor is $(22,29)$, is
\[\begin{array}{rl}
\M(\Lambda_{\{1,3\}})=&\{(3,4),(6,8),(7,9),(9,12),(10,14),(12,16),(13,19),(14,18),\\
                      & (15,20),(17,23),(18,24),(21,28)\}.
											\end{array}
\] 
Notice that $\#\M(\Lambda_{\{1,3\}})=12$.\medskip

The set of maximal points of $\Lambda_{\{2,3\}}$, whose conductor is $(22,29)$, is
\[\begin{array}{rl}
\M(\Lambda_{\{2,3\}})=&\{(3,4),(6,8),(7,9),(9,12),(10,14),(11,13),(12,16),(13,19),\\
                      & (14,18),(15,20),(17,23),(18,24),(21,28)\}.
\end{array}
\]
Notice that $\#\M(\Lambda_{\{2,3\}})=13$.\medskip

We also have that
\[
\begin{array}{rl}
\RM(\Lambda)=&\{(14,14,17),(17,17,21),(18,18,22),(20,20,25),(21,21,26),\\
             & (22,22,27),(23,23,29),(24,24,30),(25,25,32),(26,26,33),\\
             & (27,27,34),(29,29,37),(30,30,38),(33,33,42)\},\\ \\
\AM(\Lambda)=&\{(3,3,4),(6,6,8),(7,7,9),(9,9,12),(10,10,14),\\
            & (12,12,16),(13,13,19),(14,14,18),(15,15,20),\\
            & (17,17,23),(18,18,24),(21,21,28)\}.	
									\end{array}			
						\]
Notice that $\#\RM(\Lambda)=14$ and $\#\AM(\Lambda)=12$.\medskip

Moreover, $\tau_1=12$, $\tau_2=11$, $\tau_3=21$, and
\[
\begin{array}{l}
I_1=I(f_1,f_2)+I(f_1,f_3)=22+27=49, \\
I_2=I(f_2,f_1)+I(f_2,f_3)=22+27=49,\\
I_3=I(f_3,f_1)+I(f_3,f_2)=27+27=54.
\end{array}
\]

Now, by applying Formula (\ref{last}), we get
\[
\tau(\Oh)=12+11+21+ \frac{1}{2}(49+49+54) +(10+12+13)+14-12= 157.
\]

This number is validated by using the routine implemented in Singular \cite{DGPS} to calculate the Tjurina number: \bigskip

\noindent \texttt{LIB "sing.lib";}\\	
\texttt{ring r=0,(x,y),ds;} \ \ \ \ \ \ \ \ //\texttt{local ring}\\
\texttt{poly f=(y3-x7)*(y3-3*x5y-x7-x8)*(y4-2*x5y2-4*x7y-x9+x10);}\\
\texttt{tjurina(f);}\\
\texttt{//-> 157 \ \ \ \ \ \ \ \ \ \ \ \ \ \ \ //Tjurina number at 0.}


\begin{thebibliography}{XXX}


\bibitem{D87} \textsc{Delgado de la Mata, F.}, \textit{The semigroup of values of a curve singularity with several branches}, Manuscripta Math. 59, pp 347-374 (1987).

\bibitem{DGPS}
Decker, W.; Greuel, G.-M.; Pfister, G.; Sch{\"o}nemann, H.: 
\newblock {\sc Singular} {4-1-2} --- {A} computer algebra system for polynomial computations.
\newblock {http://www.singular.uni-kl.de} (2019).


\bibitem{GH} \textsc{Guzm\'an, E.M.N}; \textsc{Hefez, A.}, \textit{On the colength of fractional ideals}, arXiv: 1907.10666v1 [math. AG] 23 Jul 2019.

\bibitem{HH} \textsc{Hefez, A.}; \textsc{Hernandes, M.E.}, \textit{Computational Methods in the Local Theory of Curves}, Publica\c c\~oes Matem\'aticas, IMPA (2001).

\bibitem{HH2} \textsc{Hefez, A.}; \textsc{Hernandes, M.E.}, \textit{Standard bases for local rings of branches and their modules of differentials}, J. Symbolic Comput. Vol. 42, pp 178-191 (2007).

\bibitem{H} \textsc{Hironaka, H.}, \textit{On the arithmetic genus and the effective genera of algebraic curves}, Mem. Kyoto 80, pp 177-195 (1957).



\bibitem{Pi} \textsc{Piene, R.},  \textit{Polar classes of singular varieties}, Ann. Sci. \'Ec. Norm. Sup\'er., Serie 4, Volume 11, No. 2, pp 247-276 (1978).


\bibitem{Po} \textsc{Pol, D.}, \textit{On the values of logarithmic residues along curves}, Ann. Inst, Fourier, Tome 68, No. 2, pp 725-766 (2018).

\bibitem{Za} \textsc{Zariski, O.}, \textit{Characterization of plane algebroid curves whose module of ditterentials has maximum torsion}, Proc. Natl. Acad. Sci. USA, Vol. 56, No. 3, pp 781-786 (1966).
\end{thebibliography}
\end{document}